\documentclass[12pt]{article}
\usepackage[dvipsnames]{xcolor}
\usepackage{amsmath}
\usepackage{setspace}
\usepackage{bm}
\usepackage{bbm}
\usepackage{amssymb}
\usepackage{indentfirst}
\usepackage[superscript]{cite}
\usepackage{geometry}
\usepackage{mathtools}
\usepackage{tocloft}

\DeclareMathOperator{\sign}{sign}
\DeclarePairedDelimiter\abs{\lvert}{\rvert}%

\newcommand*{\citena}[1]{%
\begingroup
[\color{Green}
\romannumeral-`\x 
\setcitestyle{numbers}%
\cite{#1}%
\endgroup
]\ignorespacesafterend
}

\newcommand*{\citesup}[1]{%
\begingroup
\color{Green}
\cite{#1}%
\endgroup
\ignorespacesafterend
}

\newcommand*{\eqrefe}[1]{%
\begingroup
(\color{BrickRed}
\romannumeral-`\x 
\setcitestyle{numbers}%
\ref{eq:#1}%
\endgroup
)\ignorespacesafterend
}

\newcommand*{\secrefe}[1]{%
\begingroup
(\color{Aquamarine}
\romannumeral-`\x 
\setcitestyle{numbers}%
\ref{#1}%
\endgroup
)\ignorespacesafterend
}

\newtagform{Tags}[\textcolor{BrickRed}]{\color{Black}(}{)}

\usepackage[super]{natbib}
\setcitestyle{super}

\newcommand{\ii}{\bm{i}}

\begin{document}
\title{The Lerch $\Phi$ Analytic Continuation}
\date{January 12, 2021}
\author{Jose Risomar Sousa}
\maketitle
\usetagform{Tags}

\begin{abstract}
We demonstrate how to extend formulae for the Lerch transcendent function, $\Phi(e^z,k,b)$, and the polylogarithm, $\mathrm{Li}_{k}(e^{z})$, that only hold at the positive integers to the right half of the complex $k$-plane, that is, $\Re{(k)}>0$. The same is done for the partial sums of each of these functions.
\end{abstract}

\tableofcontents

\section{Introduction}
Formulae for the Lerch $\Phi$ and the polylogarithm functions were produced in \citena{Lerch}, which are only valid for positive integer $k$. Now they are extended to the right half of the complex $k$-plane, that is, $\Re{(k)}>0$.\\

These formulae provide the analytic continuation of the infinite series,
\begin{equation} \label{eq:inf_series}
\Phi(e^z,k,b)=\sum _{j=0}^{\infty}\frac{e^{z\,j}}{(j+b)^{k}} \text{, and } \mathrm{Li}_{k}(e^{z})=\sum _{j=1}^{\infty}\frac{e^{z\,j}}{j^{k}} \text{,}
\end{equation}
\noindent to the complex $z$-plane, since they hold even when the real part of $z$ is greater than zero (as long as the imaginary part of $z$ is not greater than $2\pi$ in modulus).\\

An expression for the Lerch $\Phi$ series can be created using the Abel-Plana formula, but it does not hold in the whole complex plane, unless one of its integrals is analytically continued, which is not trivial to do.\\

The new formulae hereby presented have a possible advantage over the Abel-Plana formula, as they express the infinite series by means of a closed-form and an integral, which is arguably more interesting than the Abel-Plana integral.\\

Throughout this text, when $k$ is not a non-negative integer, $k!$ is to be understood as $\Gamma{(k+1)}$, which is defined for all complex $k$, except the negative integers.

\section{The formulae at the positive integers}
The starting point for the formula derivations that follow next are the Lerch's transcendent, $\Phi(e^z,k,b)$, and the polylogarithm, $\mathrm{Li}_{k}(e^{z})$, formulae created in \citena{Lerch}, which are presented in the subsequent sections. They hold for all integer $k\ge 1$ and all complex $z$ and $b$, unless otherwise noted.\\

The partial and full polylogarithm formulae stem from the Lerch $\Phi$ formulae. Both are derived by taking the limit of the respective $\Phi$ formulae as $b$ tends to zero.

\subsection{Partial Lerch $\Phi$}
The partial Lerch's $\Phi$ sums are given by,
\begin{multline} \label{eq:lerch_parc} 
\sum _{j=1}^{n}\frac{e^{z(j+b)}}{(j+b)^{k}}=-\frac{e^{z\,b}}{2b^{k}}+\frac{e^{z(n+b)}}{2(n+b)^{k}}+\frac{1}{2b^{k}}\sum_{j=0}^{k}\frac{(z\,b)^j}{j!}-\frac{1}{2(n+b)^{k}}\sum_{j=0}^{k}\frac{(z (n+b))^j}{j!}\\+\sum_{j=1}^{k}\frac{z^{k-j}}{(k-j)!}HP_{j}(n)+\frac{z^{k}}{2(k-1)!}\int_0^1(1-u)^{k-1}\left(e^{z(n+b) u}-e^{z\,b\,u}\right)\coth{\frac{z u}{2}}\,du \text{,}
\end{multline}
\noindent where the $HP_{j}(n)$ are the so-called generalized harmonic progressions,
\begin{equation} \nonumber
HP_{j}(n)=\sum_{q=1}^{n}\frac{1}{(q+b)^j} \text{}
\end{equation}

\subsection{Full Lerch $\Phi$}
The Lerch $\Phi$ function is the limit of the previous expression when $n$ tends to infinity. The below limits hold for all complex $z$, except $z$ with non-negative real part and absolute imaginary part greater than $2\pi$ (that is, $\Re{(z)}>=0$ and $\abs{\Im{(z)}}>2\pi$).\\

Since the main expression is improper for integer and half-integer $b$, alternative expressions are provided. Note there is a slight transformation in relation to the end results shown in \citena{Lerch}.

\subsubsection{Non-integer $2b$}
For all integer $k \ge 1$,
\begin{multline} \label{eq:lerch_full} 
\sum _{j=1}^{\infty}\frac{e^{z(j+b)}}{(j+b)^{k}}=-\frac{1}{2b^{k}}\left(e^{z\,b}+\sum_{j=0}^{k-2}\frac{(z\,b)^j}{j!}\right)+\sum_{j=2}^{k}\frac{z^{k-j}}{(k-j)!}\zeta(j,b)\\+\frac{\pi\,z^k}{2(k-1)!}\cot{\pi b}-\frac{z^{k-1}}{(k-1)!}\log{\left(-\frac{z}{2\pi}\right)}\\-\frac{z^{k}}{2(k-1)!}\int_0^1(1-u)^{k-1}e^{z\,b\,u}\coth{\frac{z u}{2}}+\frac{2\pi}{z}\left(-1+\frac{\sin{2\pi b u}}{\sin{2\pi b}}\right)\cot{\pi u}\,du 
\end{multline}

\subsubsection{Half-integer $b$}
For all integer $k \ge 1$,
\begin{multline} \nonumber 
\sum _{j=1}^{\infty}\frac{e^{z(j+b)}}{(j+b)^{k}}=-\frac{1}{2b^{k}}\left(e^{z\,b}+\sum_{j=0}^{k-2}\frac{(z\,b)^j}{j!}\right)+\sum_{j=2}^{k}\frac{z^{k-j}}{(k-j)!}\zeta(j,b)\\-\frac{z^{k-1}}{(k-1)!}\log{\left(-\frac{z}{\pi}\right)}-\frac{z^{k}}{2(k-1)!}\int_0^1(1-u)^{k-1}e^{z\,b\,u}\coth{\frac{z u}{2}}-\frac{\pi}{z}\cos{\pi b u}\cot{\frac{\pi u}{2}}\,du 
\end{multline}

\subsubsection{Integer $b$}
For all integer $k \ge 1$,
\begin{multline} \nonumber 
\sum _{j=1}^{\infty}\frac{e^{z(j+b)}}{(j+b)^{k}}=-\frac{1}{2b^{k}}\left(e^{z\,b}+\sum_{j=0}^{k-2}\frac{(z\,b)^j}{j!}\right)+\sum_{j=2}^{k}\frac{z^{k-j}}{(k-j)!}\zeta(j,b)\\-\frac{z^{k-1}}{(k-1)!}\log{\left(-\frac{z}{2\pi}\right)}-\frac{z^{k-1}}{(k-1)!}\left(H(b)-\frac{1}{2b}\right)\\-\frac{z^{k}}{2(k-1)!}\int_0^1(1-u)^{k-1}e^{z\,b\,u}\coth{\frac{z u}{2}}-\frac{2\pi}{z}(1-u)\cot{\pi u}\,du 
\end{multline}

\subsection{Partial polylogarithm} \label{part2}
If one takes the limit of equation \eqrefe{lerch_parc} as $b$ tends to 0, one  obtains,
\begin{multline} \label{eq:polylog_parc}
\sum _{j=1}^n \frac{e^{z\,j}}{j^{k}}=\frac{e^{z\,n}}{2n^{k}}-\frac{1}{2n^{k}}\sum_{j=0}^{k}\frac{(z\,n)^j}{j!}+\sum_{j=1}^{k}\frac{z^{k-j}}{(k-j)!}H_{j}(n)\\
+\frac{z^{k}}{2(k-1)!}\int_0^1(1-u)^{k-1}\left(e^{z\,n u}-1\right)\coth{\frac{z u}{2}}\,du 
\end{multline}
\indent This holds for all complex $z$. Let us call these partial polylog sums $E^z_k(n)$.

\subsection{Full polylogarithm} \label{full_poly}
The following expression for the polylogarithm is obtained taking the limit of \eqrefe{lerch_full} as $b$ tends to 0. It holds for all integer $k\ge 1$ and all complex $z$, except $z$ such that $\Re{(z)}>=0$ and $\abs{\Im{(z)}}>2\pi$,
\begin{multline} \label{eq:polylog_full}
\sum _{j=1}^{\infty} \frac{e^{z\,j}}{j^{k}}=-\frac{z^{k}}{2k!}-\frac{z^{k-1}}{(k-1)!}\log{\left(-\frac{z}{2\pi}\right)}+\sum_{j=2}^{k}\frac{z^{k-j}}{(k-j)!}\zeta(j)\\
-\frac{z^{k}}{2(k-1)!}\int_{0}^{1}(1-u)^{k-1}\coth{\frac{z u}{2}}-\frac{2\pi}{z}(1-u)\cot{\pi u}\,du
\end{multline}

\section{Polylogarithm extension} \label{polylog}
The demonstration starts with the polylogarithm formulae and then a generalization for the Lerch $\Phi$ function follows.

\subsection{Partial polylogarithm} \label{AC_part}
As incredible as it may seem, if the two finite sums in the formula \eqrefe{polylog_parc} are replaced with closed-forms (more precisely, integrals), the result happens to be its extension. However, the way that is done affects the domain of the extended formula and the approach presented here seems to be ideal.\\

First we have the finite sum,
\begin{equation} \nonumber
\sum_{j=1}^{k}\frac{z^{k-j}}{(k-j)!}H_{j}(n)=z^k\sum_{j=1}^{k}\frac{z^{-j}}{(k-j)!}\sum_{q=1}^{n}\frac{1}{q^j}=z^k\sum_{q=1}^{n}\sum_{j=1}^{k}\frac{(z\,q)^{-j}}{(k-j)!}
\end{equation}\\
\indent Let us take just the part that matters,
\begin{equation} \nonumber
\sum_{j=1}^{k}\frac{(z\,q)^{-j}}{(k-j)!}=\frac{-1+e^{z\,q}(z\,q)^{-k}\,\Gamma{(k+1,z\,q)}}{k!} \text{,}
\end{equation}
\noindent where $\Gamma$ is the incomplete gamma function, which is given by an integral.\\

Let us then rewrite the initial sum,
\begin{equation} \nonumber
\sum_{q=1}^{n}\sum_{j=1}^{k}\frac{(z\,q)^{-j}}{(k-j)!}=\sum_{q=1}^{n}\left(-\frac{1}{k!}+\frac{e^{z\,q}}{k!\,(z\,q)^{k}}\int_{z\,q}^{\infty}t^k e^{-t}\,dt\right) \text{}
\end{equation}\\
\indent This last integral needs to be transformed as below so it can be summed over $q$,
\begin{equation} \nonumber
\Gamma{(k+1,z\,q)}=\int_{z\,q}^{\infty}t^k e^{-t}\,dt=z\,q\int_{1}^{c\,\infty}(z\,q\,t)^k e^{-z\,q\,t}\,dt=\int_{0}^{1}(z\,q-\log{u})^k\,e^{-z\,q}\,du \text{,}
\end{equation}
\noindent where $c$ is the directional angle of the transformation (since $z$ is complex) given by,
\begin{equation} \nonumber
c=\begin{cases}
\sign{\Re{(z)}}, & \text{if }\Re{(z)} \neq 0\\
-\ii\sign{\Im{(z)}}, & \text{otherwise.}
\end{cases}
\end{equation}
\indent Therefore, introducing the constant to the integral upper endpoint,
\begin{equation} \nonumber
\int_{1}^{c\,\infty}\sum_{q=1}^{n}\left(-\frac{z\,q\,e^{z\,q(1-t)}}{k!}+\frac{z\,q\,t^k\,e^{z\,q(1-t)}}{k!}\right)\,dt \text{,}
\end{equation}
\noindent and collapsing the sum over $q$ it can be concluded that,
\begin{equation} \nonumber
\sum_{j=1}^{k}\frac{z^{-j}}{(k-j)!}H_j(n)=\frac{z}{k!}\int_{1}^{c\,\infty}(t^k-1)\frac{e^{z(1-t)}+n\,e^{z(1-t)(n+2)}-(n+1)e^{z(1-t)(n+1)}}{(e^{z(1-t)}-1)^2}\,dt
\end{equation}
\indent Now the integral is transformed to make the constant $c$ disappear,
\begin{equation} \nonumber
\sum_{j=1}^{k}\frac{z^{k-j}}{(k-j)!}H_j(n)=\frac{z^k}{k!}\int_{0}^{1}\frac{1-(n+1)u^n+n\,u^{n+1}}{(1-u)^2}\left(-1+z^{-k}\left(z-\log{u}\right)^k\right)\,du
\end{equation}
\indent Note that due to to the presence of multi-valued functions such as $z^{-k}$, when $k$ is for instance a real number, it is not advisable to simplify the expression within the above integral.\\

The second and last finite sum for which a closed-form is needed is given by,
\begin{equation} \nonumber
\sum_{j=0}^{k}\frac{(z\,n)^j}{j!}=\frac{e^{z\,n}\,\Gamma{(k+1,z\,n)}}{k!}=\frac{1}{k!}\int_{0}^{1}\left(z\,n-\log{u}\right)^k\,du \text{,}
\end{equation}
\noindent which concludes the exercise.\\

For all complex $k$ and $z$ such that $\Re{(k)}>0$,
\begin{multline} \label{eq:E^z_k(n)}
\sum _{j=1}^n \frac{e^{z\,j}}{j^{k}}=\frac{e^{z\,n}}{2n^{k}}-\frac{e^{z\,n}}{2n^k}\frac{\Gamma{(k+1,z\,n)}}{k!}\\
-\frac{n\,z^k}{k!}+\frac{1}{k!}\int_{0}^{1}\frac{1-(n+1)u^n+n\,u^{n+1}}{(1-u)^2}\left(z-\log{u}\right)^k\,du\\
+\frac{z^{k}}{2(k-1)!}\int_0^1(1-u)^{k-1}\left(e^{z\,n\,u}-1\right)\coth{\frac{z u}{2}}\,du \text{}
\end{multline}
\indent The first term on the second line came from a simplification of the former integral.

\subsection{$H_k(n)$ when $\Re{(k)}>-1$} 
From formula \eqrefe{E^z_k(n)}, different formulae can be derived for the generalized harmonic numbers\citesup{GHN}, $H_k(n)$. The formula allows one degree of freedom ($z$ can be any multiple of $2\pi\ii$). When $z=0$, the simplest possible expression is obtained, which surprisingly holds for $\Re{(k)}>-1$,
\begin{equation} \label{eq:H_k(n)}
\sum_{j=1}^n\frac{1}{j^{k}}=\frac{1}{k!}\int_{0}^{1}\left(\frac{1-(n+1)u^n+n\,u^{n+1}}{(1-u)^2}\right)(-\log{u})^k\,du
\end{equation}
\indent The limit of \eqrefe{H_k(n)} when $n$ approaches infinity is $\zeta(k)$, whose expression is given in \eqrefe{zeta(k)}. Hence, if $\Re{(k)}>1$, $H_k(n)$ can be rewritten as,
\begin{equation} \nonumber
\sum_{j=1}^n\frac{1}{j^{k}}=\zeta{(k)}-\frac{1}{k!}\int_{0}^{1}\frac{(n+1)u^n-n\,u^{n+1}}{(1-u)^2}(-\log{u})^k\,du
\end{equation}

\subsection{$\zeta(k)$ outside the critical strip} \label{zeta}
The limit of \eqrefe{H_k(n)} when $n$ goes to infinity is $\zeta(k)$. Another way to derive this formula is to sum the delta operator, $\Delta{H_k(n)}$, applied to \eqrefe{H_k(n)}, over all the naturals. If $\Re{(k)}>1$,
\begin{equation} \label{eq:zeta(k)}
\zeta(k)=\frac{1}{k!}\int_{0}^{1}\frac{\left(-\log{u}\right)^k}{(1-u)^2}\,du=\sum_{n=1}^{\infty}\Delta{H_k(n)}=\frac{1}{(k-1)!}\int_{0}^{1}\frac{\left(-\log{u}\right)^{k-1}}{1-u}\,du 
\end{equation}
\indent The rightmost formula coincides with the one from the literature.\\

If $\Re{(k)}<0$, the zeta function has a different integral representation, obtained via Riemann's functional equation,
\begin{equation} \nonumber
\zeta(k)=-\frac{2(2\pi)^{k-1}}{k-1}\sin{\frac{k\,\pi}{2}}\int_{0}^{1}\frac{\left(-\log{u}\right)^{-k+1}}{(1-u)^2}\,du 
\end{equation}

\subsection{Full polylogarithm}
For the purpose of producing the extension of $\mathrm{Li}_{k}(e^{z})$, it is more convenient to start from its formula seen in \secrefe{full_poly} and repeat some of the previous steps from \secrefe{AC_part}. For the polylogarithm, the first and only finite sum is,
\begin{equation} \nonumber
\sum_{j=2}^{k}\frac{z^{k-j}}{(k-j)!}\zeta(j)=z^k\sum_{q=1}^{\infty}\sum_{j=2}^{k}\frac{(z\,q)^{-j}}{(k-j)!} \text{,}
\end{equation}
\noindent and to find out the closed-form one can take the limit of the partial sum as $n$ goes to infinity,
\begin{equation} \nonumber
\sum_{j=2}^{k}\frac{z^{-j}}{(k-j)!}\zeta(j)=\lim_{n\to\infty}\int_{1}^{c\,\infty}\sum_{q=1}^{n}\left(-\frac{e^{z\,q(1-t)}}{(k-1)!}+\frac{z}{k!}\left(t^k-1\right)q\,e^{z\,q(1-t)}\right)\,dt \text{}
\end{equation}\\
\indent Next, the sum inside the integral can be simplified collapsing the sum over $q$,
\begin{equation} \nonumber
\lim_{n\to\infty}\int_{1}^{c\,\infty}\left(-\frac{e^{z(1-t)}-e^{z(1-t)(n+1)}}{(k-1)!\left(e^{z(1-t)}-1\right)}+z\left(t^k-1\right)\frac{e^{z(1-t)}+n\,e^{z(1-t)(n+2)}-(n+1)e^{z(1-t)(n+1)}}{k!\left(e^{z(1-t)}-1\right)^2}\right)\,dt
\end{equation}\\
\indent This limit is not trivial. The solution might be found if one assumes that all terms with $n$ cancel out (which is a plausible assumption if, for example, $z$ is a positive real number) and the remaining terms constitute the limit.\\

It turns out this hunch is only half right, but fortunately it is possible to find out the term that makes up the difference through trial and error, leading to the below,
\begin{equation} \nonumber
\sum_{j=2}^{k}\frac{z^{-j}}{(k-j)!}\zeta(j)=-\frac{1}{z(k-1)!}+\frac{e^z}{k!}\int_{1}^{c\,\infty}\frac{k\,e^{z}+(-k-z+z\,t^k)e^{z\,t}}{\left(e^{z}-e^{z\,t}\right)^2}\,dt
\end{equation}
\indent It is best to make a change of variables ($u=e^{-z(t-1)}$) in order to do away with $c$ and enable the combination of the integrals,
\begin{equation} \nonumber
\sum_{j=2}^{k}\frac{z^{k-j}}{(k-j)!}\zeta(j)=-\frac{z^{k-1}}{(k-1)!}+\frac{z^{k-1}}{k!}\int_{0}^{1}\frac{-k\,u-z+z^{-k+1}\left(z-\log{(1-u)}\right)^k}{u^2}\,du
\end{equation}\\
\indent Therefore, when all is put together one finds that,
\begin{multline} \nonumber
\mathrm{Li}_{k}(e^{z})=-\frac{z^{k}}{2k!}-\frac{z^{k-1}}{(k-1)!}\log{\left(-\frac{e\,z}{2\pi}\right)}+\\
-\frac{z^{k}}{2(k-1)!}\int_{0}^{1}(1-u)^{k-1}\coth{\frac{z u}{2}}-\frac{2\pi}{z}(1-u)\cot{\pi u}+\frac{2(k\,u+z-z(1-\log{(1-u)/z})^k}{k\,z\,u^2}\,du
\end{multline}\\
\indent It is possible to simplify the above integral with the observation that,
\begin{equation} \nonumber
\int_{0}^{1}-\frac{2\pi}{z}(1-u)\cot{\pi u}+\frac{2}{z\,u}\,du=\frac{2\log2\pi}{z} \text{}
\end{equation}\\
\indent Finally, for all $k$ such that $\Re{(k)}>0$ and all $z$ (except $\Re{(z)}>=0$ and $\abs{\Im{(z)}}>2\pi$),
\begin{multline} \label{eq:Polylog_final}
\mathrm{Li}_{k}(e^{z})=-\frac{z^{k}}{2k!}-\frac{z^{k-1}\left(1+\log{(-z)}\right)}{(k-1)!}\\
-\frac{z^{k}}{2(k-1)!}\int_{0}^{1}(1-u)^{k-1}\coth{\frac{z u}{2}}+\frac{2}{k\,u^2}\left(1-z^{-k}\,(z-\log(1-u))^k\right)\,du
\end{multline}
\indent It is important to note that if $\abs{\Im{(z)}}=2\pi$ one must have $\Re{(k)}>1$. This convergence domain has been thoroughly checked but might still be subject to change.

\section{Lerch $\Phi$ extension}
There is a shortcut to obtain these formulae, which uses a relation between the Lerch $\Phi$ and the full and the partial polylogs that is trivial to demonstrate. Here the detailed steps provided in section \secrefe{polylog} are omitted, since they can be easily replicated by analogy, but some striking differences between the two cases are pointed out.

\subsection{Partial Lerch $\Phi$ series}
In the case of the partial Lerch $\Phi$ series, following the same process employed for the polylogarithm, we conclude that the sum of $HP_j(n)$ can be extended by,
\begin{equation} \label{eq:hp_soma} \nonumber
\sum_{j=1}^{k}\frac{z^{k-j}}{(k-j)!}HP_j(n)=\frac{z^k}{k!}\int_{0}^{1}u^b\left(b\frac{1-u^n}{1-u}+\frac{1-(n+1)u^n+n\,u^{n+1}}{(1-u)^2}\right)\left(-1+z^{-k}\left(z-\log{u}\right)^k\right)\,du
\end{equation}
\indent The above integral converges when $\Re{(b)}>-1$, if $k$ is a non-negative integer. The complement of this domain would require a much more convoluted expression, so we are not trying to solve it. This is justified by the fact that the constant $c$ from section \secrefe{AC_part} now depends on the summation index $q$ (e.g., $\sign{\Re{(z(q+b))}}$). If, however, $k$ is not a non-negative integer, then one must have $\Re{(k)}>0$ and $\Re{(b)}>0$ for the integral to converge.\\ 

The integral can be further simplified as,
\begin{equation} \nonumber
\sum_{j=1}^{k}\frac{z^{k-j}}{(k-j)!}HP_j(n)=-\frac{n\,z^k}{k!}+\frac{1}{k!}\int_{0}^{1}u^b\left(b\frac{1-u^n}{1-u}+\frac{1-(n+1)u^n+n\,u^{n+1}}{(1-u)^2}\right)\left(z-\log{u}\right)^k\,du
\end{equation}\\
\indent Therefore, the partial Lerch $\Phi$ function for complex $k$, $z$ and $b$, such that $\Re{(k)}>0$ and $\Re{(b)}>0$ ($\Re{(b)}>-1$, if $k$ is integer), is given by,
\begin{multline} \label{eq:partial_Lerch_AC}
\sum _{j=1}^{n}\frac{e^{z(j+b)}}{(j+b)^{k}}=-\frac{e^{z\,b}}{2b^{k}}+\frac{e^{z\,b}}{2b^k}\frac{\Gamma{(k+1,z\,b)}}{k!}+\frac{e^{z(n+b)}}{2(n+b)^{k}}-\frac{e^{z(n+b)}}{2(n+b)^k}\frac{\Gamma{\left(k+1,z(n+b)\right)}}{k!}\\-\frac{n\,z^k}{k!}+\frac{1}{k!}\int_{0}^{1}u^b\left(b\frac{1-u^n}{1-u}+\frac{1-(n+1)u^n+n\,u^{n+1}}{(1-u)^2}\right)\left(z-\log{u}\right)^k\,du\\+\frac{z^{k}}{2(k-1)!}\int_0^1(1-u)^{k-1}\left(e^{z(n+b) u}-e^{z\,b\,u}\right)\coth{\frac{z u}{2}}\,du \text{}
\end{multline}

\subsection{$HP_{k}(n)$ when $\Re{(k)}>-1$}
As a consequence of equation \eqrefe{partial_Lerch_AC}, an alternative integral representation for the generalized harmonic progressions, $HP_k(n)$\citesup{GHP}, valid when $\Re{(k)}>-1$ and $\Re{(b)}>-1$, is,
\begin{equation} \label{eq:HP_k(n)}
\sum_{j=1}^{n}\frac{1}{(j+b)^k}=\frac{1}{k!}\int_{0}^{1}u^b\left(b\frac{1-u^n}{1-u}+\frac{1-(n+1)u^n+n\,u^{n+1}}{(1-u)^2}\right)\left(-\log{u}\right)^k\,du
\end{equation}\\
\indent Now, using the Hurwitz zeta formula from \eqrefe{Hurwitz}, if $\Re{(k)}>1$ and $\Re{(b)}>0$, $HP_k(n)$ can be rewritten as,
\begin{equation} \nonumber
\sum_{j=1}^{n}\frac{1}{(j+b)^k}= \zeta(k,b+1)-\frac{1}{k!}\int_{0}^{1}u^b\left(\frac{b\,u^n}{1-u}+\frac{(n+1)u^n-n\,u^{n+1}}{(1-u)^2}\right)\left(-\log{u}\right)^k\,du
\end{equation}

\subsection{$\zeta(k,b+1)$ when $\Re{(k)}>1$}
From the limit of \eqrefe{HP_k(n)} as $n$ tends to infinity, the Hurwitz zeta function is obtained. If $\Re{(k)}>1$ and $\Re{(b)}>0$ ($\Re{(b)}>-1$ in some cases),
\begin{equation} \label{eq:Hurwitz}
\zeta(k,b+1)=\frac{1}{k!}\int_{0}^{1}u^b\left(\frac{b}{1-u}+\frac{1}{(1-u)^2}\right)\left(-\log{u}\right)^k\,du
\end{equation}

\subsection{Full Lerch $\Phi$ series}
For the full Lerch $\Phi$ series, the sum of the Hurwitz zeta functions needs to extended to values of $k$ other than the positive integers. Like before, for  non-negative integer $k$ and every $z$ the below formula holds when $\Re{(b)}>-2$,
\begin{multline} \label{eq:hurwitz_soma} \nonumber
\sum_{j=2}^{k}\frac{z^{k-j}}{(k-j)!}\zeta(j,b)=-\frac{z^k}{k!}-\frac{z^{k-1}}{(k-1)!}\left(1+\frac{1}{b}\right)+\frac{e^{z\,b}}{b^k}\frac{\Gamma{\left(k+1,z\,b\right)}}{k!}\\-\frac{z^{k-1}}{k!}\int_{0}^{1}\frac{(1-u)^{b}}{u^2}\left(k\,u+z\left(1+b\,u\right)\left(1-z^{-k}\left(z-\log(1-u)\right)^k\right)\right)\,du
\end{multline}
\indent This sum could have been obtained more easily using equation \eqrefe{Hurwitz} as a shortcut.\\

When the above expression is replaced in equation \eqrefe{lerch_full}, the resulting formula can be simplified using the below equation, valid if $\Re{(b)}>-1$,
\begin{equation} \nonumber 
\int_0^1 \pi\left(-1+\frac{\sin{2\pi b u}}{\sin{2\pi b}}\right)\cot{\pi u}+\frac{(1-u)^{b}}{u}\,du=\log{2\pi}-\frac{1}{2b}+\frac{\pi\cot{\pi b}}{2} \text{}
\end{equation}
\indent This simplification is convenient, since it does away with $\sin{2\pi\,b}$ in the fraction's denominator within the integral, which is improper for integer and half-integer $b$, and also since the trigonometric parts of the formula cancel out.\\

Therefore, the below formula for the Lerch transcendent should be valid when $\Re{(k)}>0$, and $z$ is not such that $\Re{(z)}\ge 0$ and $\abs{\Im{(z)}}>2\pi$, and $\Re{(b)}>0$,
\begin{multline} \nonumber \label{eq:Lerch_final} 
e^{z\,b}\,\Phi(e^z,k,b)=-\frac{z^k}{2\,k!}-\frac{z^{k-1}\left(1+\log{(-z)}\right)}{(k-1)!}+\frac{e^{z\,b}}{2\,b^k}\left(1+\frac{\Gamma{\left(k+1,z\,b\right)}}{k!}\right)
\\-\frac{z^{k}}{2(k-1)!}\int_0^1(1-u)^{k-1}e^{z\,b\,u}\coth{\frac{z u}{2}}+\frac{2(1-u)^{b}\left(1+b\,u\right)}{k\,u^2}\left(1-z^{-k}\left(z-\log(1-u)\right)^k\right)\,du 
\end{multline}
\indent The restriction on $b$, $\Re{(b)}>0$, did not exist for integer $k$.


\end{document}